\input amstex
\input epsf.sty
\documentstyle{amsppt}
\document
\magnification=1200
\NoBlackBoxes
\nologo
\pageheight{18cm}

%\hfill{\it file Amywork/gonchnew.tex, 25.04.02}

\bigskip

\centerline{\bf MULTIPLE $\zeta$--MOTIVES AND MODULI SPACES 
$\overline{\Cal{M}}_{0,n}$}

\bigskip

\centerline{A.~B.~Goncharov${}^{(1)}$, Yu.~I.~Manin${}^{(2)}$}

\medskip

\centerline{\it ${}^{(1)}$ Brown University, Providence RI, USA}

\smallskip

\centerline{\it ${}^{(2)}$ Max--Planck--Institut f\"ur Mathematik, Bonn, Germany}

\bigskip

{\bf Contents}

\medskip

1\ Introduction and summary

2\ Geometry of $\overline{\Cal{M}}_{0,n+3}$ and divisors $B$

3\ Forms  $\Omega(\underline \varepsilon )$
and divisors $A(\underline \varepsilon )$ 

4\ Multiple $\zeta$--motives

\bigskip

{\bf 1\quad Introduction and summary}

\medskip

{\bf 1. $\zeta$--values and moduli spaces.} (i) The
multiple $\zeta$-values were defined by Euler [E] as the series
$$
\zeta(n_1, ..., n_m) = \sum_{0< k_1 < ... < k_m}
\frac{1}{k_1^{n_1} ...\, k_m^{n_m} }\, , \qquad n_m >1 \,.
\eqno(1)
$$

(ii) Recall  that $\overline{{\Cal M}}_{0, n+3}$ parametrizes
stable curves of genus 0 with $n+3$ labeled points.
It is stratified according to the type of degeneration.
The open stratum ${\Cal M}_{0, n+3}$
can be identified  with
$(\bold{P}^1)^{n+3}$ minus diagonals modulo $\roman{Aut}\,\bold{P}^1$.
Then
$$
\partial \overline {\Cal M}_{0, n+3}:= \overline {\Cal M}_{0, n+3} -
{\Cal M}_{0, n+3}
$$
is a normal crossing divisor, and the pair
 $(\overline {\Cal M}_{0, n+3}, \partial \overline
{\Cal M}_{0, n+3})$ is defined
over $\bold Z$. We use it to construct an unramified over $\bold Z$   framed
mixed Tate motive whose period is given by the
multiple $\zeta$--value (1).

\smallskip

Generally, let $A$ and $B$ be unions of some irreducible components
of $\partial \overline {\Cal M}_{0, n+3}$. Suppose that
no irreducible component is shared by   $A$ and $B$.
Then we show that
$$
H^n(\overline {\Cal M}_{0, n+3} -A, B - A \cap B)
\eqno(2)
$$
is   a mixed Tate motive  unramified
over $\bold Z$.  A choice of  two non zero elements
$$
[\Omega_A] \in \roman{Gr}^W_{2n}H^n( \overline {\Cal M}_{0, n+3} -A);
\quad
[\Delta_B] \in \left(\roman{Gr}^W_{0}H^n(\overline {\Cal M}_{0,
n+3},B)\right)^{\vee}
\eqno(3)
$$
defines {\it a framed} mixed Tate motive unramified over $\bold Z$ 
given by the triple
$$
\left(H^n(\overline {\Cal M}_{0, n+3} -A, B_A); [\Omega_A],
[\Delta_B]\right) \, .
\eqno(4)
$$

For every multiple $\zeta$--value (1) we  construct
two divisors  on
$\overline {\Cal M}_{0, n+3}$  as above such  that the  period
of the corresponding motive is this value.

\smallskip

Another
construction of the multiple polylogarithm motives has been
given in chapter 3 of [G5]. The construction given here differs
from the earlier one for multiple $\zeta$--values.

\medskip

{\bf 2. A sketch of the construction of multiple $\zeta$-motives}.
First, let us recall how to represent multiple $\zeta$--values by
iterated integrals.
Starting with positive integers $n_1,\dots ,n_m$ as in (1), we put
$n:=n_1+\dots +n_m$, and  $\underline \varepsilon:=
(\varepsilon_1, ..., \varepsilon_n)$ where $\varepsilon_i = 0$ or $1$,
and
$\varepsilon_i=1$ precisely when
$i\in\{1,n_1+1,n_1+n_2+1,\dots ,n_1+\dots + n_{m-1}+1\}$.
Furthermore, put
$$
\omega(\underline \varepsilon):= \frac{dt_1}{t_1-\varepsilon_1}\wedge
...
\wedge \frac{dt_n}{t_n-\varepsilon_n}
$$
and
$$
\Delta^0_n:= \{(t_1^0, \dots, t_n^0) \in \bold{R}^n\quad | \quad 0 < t_1^0 < \dots <
t_n^0 < 1\}
$$
By the Leibniz-Kontsevich formula (cf. Theorem 2.2 in [G4])
$$
\zeta (n_1,\dots ,n_m)= \zeta (\underline \varepsilon )
= (-1)^m\int_{\Delta^0_n} \omega (\underline \varepsilon ).
$$
The integral converges because $\varepsilon_1 \not = 0, \varepsilon_n
\not = 1$.

\medskip

Let us identify points of ${\Cal M}_{0,n+3}(\bold{C})$
with sequences $(t_1^0,\dots ,t_n^0)\in\bold{C}^n$ such that $t_i^0\ne 0,1$,
$t_i^0\ne t_j^0$ for $i\ne j$, using the rule
$$
\{0, t_1^0,\dots ,t_n^0, 1, \infty\}\ \roman{mod\ Aut}\,\bold{P}^1\ \Leftrightarrow\ (t_1^0,\dots ,t_n^0)
$$
Moreover, consider the
form $\omega (\underline{\varepsilon})$ as
a holomorphic form on ${\Cal M}_{0, n+3}$. Meromorphically extending it to
$\overline {\Cal M}_{0,n+3}$ we get  a form
$\Omega(\underline{\varepsilon})$
with logarithmic singularities at the boundary
$\partial \overline {\Cal M}_{0,n+3}$.
Finally, we identify the open simplex $\Delta^0_n$ with its image
$\Phi(\Delta^0_n)$
in ${\Cal M}_{0, n+3}$. Then we obviously have
$$
\zeta (n_1,\dots ,n_m) = (-1)^m\int_{\Phi(\Delta^0_n)}\Omega (\underline
\varepsilon )
\eqno(5)
$$
Now let us introduce divisors $A = A(\underline \varepsilon )$ and $B$
and
interpret this integral as a period of a framed mixed Hodge
structure (cf. section 4 of [G3] and Chapter 4).
\smallskip

{\it The divisor $A(\underline \varepsilon)$}. It
 is the divisor of singularities of the form $\Omega (\underline
\varepsilon)$.

\smallskip

{\it The divisor $B$}.
Let ${\overline \Phi_n}$ be the closure of the open simplex
$\Phi(\Delta^0_n)$.  We define
$B$ as  the Zariski closure of the  boundary of $\overline \Phi_n$.
Thus there is a relative homology class
$$
[\overline \Phi_n] \in H_n(\overline {\Cal M}_{0, n+3}(\bold{C}), B(\bold{C}); \bold{Z})
$$
In the Chapters 2 and  3 we will clarify the structure of arbitrary choices
involved in these constructions. Moreover, we will
explicitly compute $B$ and $A(\underline \varepsilon)$,
in terms of the standard combinatorial description
of boundary components of $\partial \overline {\Cal M}_{0,
n+3}$.

\smallskip

The closed cell
 $\overline \Phi_n$
is naturally identified with a Stasheff polytope. The divisor $B$ is
an algebraic counterpart  of the Stasheff
polytope. Specifically, we show that 
intersecting a $k$--dimensional stratum of the divisor $B$ with
the cell $\overline \Phi_n$ we get a real $k$--dimensional
face of this cell.  This face is naturally identified with the corresponding
 face of the Stasheff polytope. This fact plus convergence of the
integral (5)
immediately imply the following key lemma.

\medskip

{\bf Lemma 1.1.}{\it
The set of the irreducible components of the divisor
$A(\underline \varepsilon)$ is disjoint with the one of $B$.
}

\medskip

In the Chapter 3
 we will directly check this Lemma using our combinatorial
calculations.

\smallskip

Lemma 1.1 guarantees that we can set
$$
\zeta^{\Cal M}(n_1, ..., n_m):= H^n\left(\overline {\Cal M}_{0, n+3}
- A(\underline \varepsilon),
B - B\,\cap\,{A(\underline \varepsilon)}\right)
$$
We show in the Chapter 4 that the form $\Omega(\underline \varepsilon)$ and
the cycle $\overline{\Phi}_n$
provide a framing, and the period of the corresponding mixed Hodge
structure is given  by (1). The fact that this period is a well defined number reflects the following
strengthening 
of Lemma 1.1:  the boundary of $\Phi_n$ does not intersect
$A(\underline \varepsilon)$ (see Corollary 3.2).

\medskip

{\bf 3. The simplest example: $\zeta^{\Cal M}(2)$}.
 The integral representation for $\zeta(2)$ is given  by the following
 formula discovered by Leibniz:
 $$
 \zeta(2)  =
 \int\int_{0<t_1 < t_2 < 1}\frac{dt_1}{1-t_1}\wedge
 \frac{dt_2}{t_2}
 $$
 Let us consider first $(t_1, t_2)$ as the affine coordinates on
  $\bold{P}^1 \times \bold{P}^1$. Then the singularities of the integrand
 provide a rectangle ${\Cal A}$ given by the punctured lines on
 the picture below,  whereas the Zariski closure of the boundary of the 
integration cycle is the bold
 triangle ${\Cal B}$ on the picture.

\smallskip

 %PICTURE

 \input epsf
 \epsfxsize = 1.5in
 \midinsert
 $$\centerline{\hbox{\epsffile{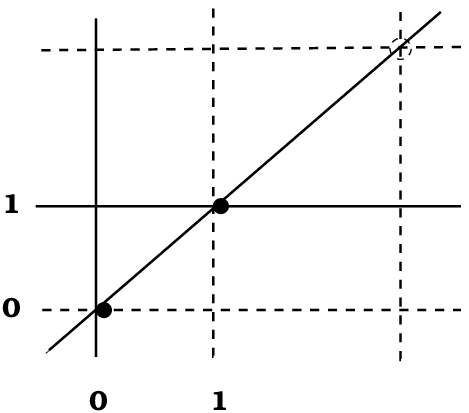}}}$$
 \endinsert

 The moduli space $\overline {\Cal M}_{0,5}$ can be obtained by blowing up
 the three special points  $(0,0), (1,1), (\infty,\infty)$ on the diagonal
 -- see the proof of the Proposition 3.1.
 They are the points where the divisor ${\Cal A} \cup {\Cal B}$
 (the union of the seven lines on the picture)
 fails to be a normal crossing divisor.

 Thus we can identify
$\overline {\Cal M}_{0,5}$ with a del Pezzo surface of degree $5$.
 The  boundary  $\partial \overline {\Cal M}_{0,5}$ is given by the ten
 exceptional curves on it. They are the strict preimages of the seven lines
 on $\bold{P}^1 \times \bold{P}^1$ plus the three  curves blowed down by
 the projection onto $\bold{P}^1 \times \bold{P}^1$.
 The divisors $A = A(1,0)$ and $B$
 are given by two complementary algebraic
 pentagons. The union of these two pentagons is the whole boundary
 $\overline {\Cal M}_{0,5}$.
 Such a decomposition of the boundary into union of two pentagons
  is not unique, but the automorphism group $S_5$ acts transitively on the
 set of all decompositions.

\smallskip

 \input epsf
 \epsfxsize = 2in
 \midinsert
 $$\centerline{\hbox{\epsffile{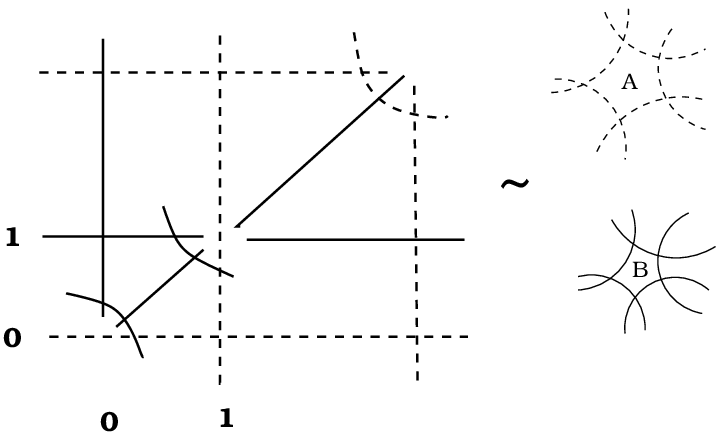}}}$$
 \endinsert

 %PICTURE

 We define the  $\zeta^{\Cal M}(2)$-motive as
 $H^2(\overline {\Cal M}_{0,5} - A, B - B \cap A)$.
 The cell $\overline \Phi_2$ is a  pentagon. It provides
 a class $[\overline \Phi_2]$ in the Betti realization of $H_2(\overline
 {\Cal M}_{0,5} - A, B - B \cap A)$. The form $\Omega(1,0)$
 gives a class in the De Rham realization of $\zeta^{\Cal M}(2)$, and the
 number
 $\zeta(2)$ equals to the pairing 
$\langle\Omega(1,0), [\overline \Phi_2]\rangle$.

\smallskip

 {\bf Remark}. It might be tempting to consider the simpler motive
 $$
H^2(\bold{P}^1\times \bold{P}^1 - {\Cal A}, {\Cal B} - {\Cal B} \cap
 {\Cal A}) \,.
$$
 However $\zeta(2)$ is not its period. Indeed,
 removing ${\Cal A}$ we delete
 a vertex  $(0,0)$ of the closed triangle $\Delta_2$. Thus
  $\Delta_2$ does not provide an element
 in the Betti realization of
 $H_2(\bold{P}^1\times \bold{P}^1 - {\Cal A}, {\Cal B} - {\Cal B} \cap
 {\Cal A})$.

\medskip

{\it Acknowledgement}.
This paper was conceived at the MPIM (Bonn) in  June of  1997,
and its rough draft was prepared during the Fall of 2001 when A.~G. visited  MSRI
(Berkeley).
A.~G. is very grateful to both
institutions for providing ideal working conditions.
The work of A.~G. was supported by the NSF grant   DMS-0099390.
We are grateful to V.~Fock  for useful discussions.

\bigskip

{\bf 2\quad Geometry of $\overline{\Cal{M}}_{0,n+3}$ and divisors $B$}

\bigskip

In this Chapter we give a summary
of the properties of $\overline{\Cal{M}}_{0,n+3}$
which will be used in the next Chapter and calculate the divisor $B$.

\medskip

{\bf 1. Reminder on moduli spaces.}
Let $S$
be any finite set, $|S|\ge 3.$
We remind that $\overline{\Cal{M}}_{0,S}$ carries
a universal family of stable curves of genus $0$,
$\pi_S:\,\overline{\Cal{C}}_{0,S}\to \overline{\Cal{M}}_{0,S}$
endowed with structure sections labeled by $S$, $x_s:\, \overline{\Cal{M}}_{0,S}\to \overline{\Cal{C}}_{0,S}$, $s\in S.$ 
Hence geometric points of $\overline{\Cal{M}}_{0,S}$ can
be identified with stable $S$--labeled pointed curves.

\smallskip

For any $S'\subset S$, $|S'|\ge 3,$ there is a canonical 
forgetful morphism
$\phi_T:\,\overline{\Cal{M}}_{0,S}\to \overline{\Cal{M}}_{0,S'}$,
$T:=S - S'.$
In terms of geometric points, $\phi_T$ forgets all sections
marked by $s\in T$ and then contracts the unstable components.

\smallskip

Morphisms forgetting disjoint subsets
of points commute in an obvious sense, and their
product is a forgetful morphism forgetting the union
of these subsets.

\smallskip

In particular, if $T=\{t\}$, we get two morphisms
with common base: $\phi_{\{t\}}$ and 
$\pi_{S'}:\,\overline{\Cal{C}}_{0,S'}\to \overline{\Cal{M}}_{0,S'}$.
A basic fact in the theory of moduli spaces
is the existence of the canonical isomorphism
$$
\mu_S:\,\overline{\Cal{C}}_{0,S'}\cong \overline{\Cal{M}}_{0,S}
\eqno(6)
$$
transforming $\phi_{\{t\}}$ to $\pi_{S'}$.

\medskip

{\bf 2. Combinatorics of boundary strata.} 
The space $\overline{\Cal{M}}_{0,S}$ is stratified by locally 
closed subschemes $\Cal{M}_{\tau}$ indexed by the (isomorphism classes of the)
stable trees with tails labeled by $S$ which are the dual graphs
of the respective stable curves. Closures of all strata are smooth.
 The open stratum
 ${\Cal{M}}_{0,S}$  parametrizes irreducible curves.
Codimension one strata parametrize two--component curves
and can be indexed by unordered 2--partitions
$S=S'\cup S''$ with $|S^{(i)}|\ge 2$ (stability), describing how the labeled points
are distributed between the two components.

\smallskip

In what follows, we will often denote such a partition as $S'|S''$
or similarly. For example, if $S=\{0,1,\infty ,s\}$, then
$0\infty |1s$ means that $S'=\{0,\infty\},\,S''=\{1,s\}.$

\smallskip

The total boundary is a normal crossing divisor. The closed
strata of the boundary are strata of this divisor:
each codimension $d$ stratum is an intersection of
exactly $d$ boundary divisors, and their set is unique.
A combinatorial description of this picture which we
will need in the following Chapter looks as follows.

\smallskip

Consider two stable 2--partitions of $S$: $\sigma :=\sigma_1|\sigma_2$
and $\tau :=\tau_1|\tau_2$. Put 
$$
\delta (\sigma ,\tau ):= \roman{the\ number\ of\ nonempty\ intersections}\
\sigma_i\,\cap\,\tau_j\ \roman{minus}\ 2.
$$
This is a kind of distance: $\delta (\sigma ,\tau )=2$ iff the respective
divisors do not intersect; 1, iff their intersection is of codimension two;
0 iff they coincide.

\smallskip

A subset $\Cal{D}$ of boundary divisors determines a non--empty stratum exactly
when all pairwise distances between the respective partitions
equal 1. Knowing this subset, we can also reconstruct the tree $\tau$
describing the respective degeneration type:
this is the maximal tree with tails labeled by $S$ such
that if one cuts it at the middle of any internal edge,
the resulting 2--partition of $S$ will belong to $\Cal{D}$.

\smallskip

Moreover, knowing the tree $\tau$ we can determine
the structure of the stratum $\Cal{M}_{\tau}$ itself: it is canonically isomorphic
to the product
$$
\prod_{v\in V_{\tau}} \Cal{M}_{0,F_{\tau}(v)}\, 
$$
where $V_{\tau}$ is the set of vertices, and $F_{\tau}(v)$
is the set of flags incident to $v$.
\smallskip

To determine the image of the boundary divisor defined by
a partition of $S$ under the
forgetful morphism $\phi_T:\,\overline{\Cal{M}}_{0,S}\to \overline{\Cal{M}}_{0,S'}$,
one should consider the induced partition of $S'$. If it
is stable, the respective divisor is the image.
If it is unstable, the image is the total space
$\overline{\Cal{M}}_{0,S'}$. For example, under the isomorphism
(6), the $S'$--labeled structure sections
of $\overline{\Cal{C}}_{0,S'}$ correspond exactly to those
boundary divisors $ts|S'-\{s\}$ of $\overline{\Cal{M}}_{0,S}$ which become
unstable after forgetting $t$.

\smallskip

For proofs and further details, see [Kn] and [Man], Ch.~3, \S 3.

\medskip

{\bf 3. Real points of $\overline{\Cal{M}}_{0,S}$.} Consider a stable $S$--labeled curve over $\bold{C}$. Endowing it with the complex conjugate
structure we will produce another similar curve. Thus
we have a conjugation involution acting upon $\overline{\Cal{M}}_{0,S}(\bold{C}).$
Fixed points of this involution form the space of
real points $\overline{\Cal{M}}_{0,S}(\bold{R}).$ A stable curve
living over a real point is itself real in the following sense:
it is endowed with a conjugation involution fixing all labeled
points and all singular points. Thus, every irreducible component
of a real stable labeled curve is a sphere $\bold{P}^1(\bold{C})$
endowed with a real equator $\bold{P}^1(\bold{R})$ carrying
all labeled points and eventually all intersection
points with other components.

\smallskip

Taken all together, such points are called {\it special ones.}
In the dual graph $\tau$, an irreducible component becomes a vertex,
and special points become flags at this vertex. Thus, a real structure
of the respective curve determines an additional structure on
$\tau$: {\it an unoriented cyclic order on the set of flags at each vertex.}
To define it, just look at how the special points are distributed
along the real equator of the respective component.

\smallskip

We will call a choice of such orders {\it a locally planar structure}
on $\tau$. 

\smallskip

Varying the curve along one connected component of the intersection
of $\overline{\Cal{M}}_{0,S}(\bold{R})$ with the $\tau$--stratum
does not change the local planar structure.

\smallskip

The following Proposition furnishes additional details
of this picture.

\medskip

{\bf Proposition 2.1.} {\it (i) The set 
$\overline{\Cal{M}}_{0,S}(\bold{R})$ is a connected closed
real manifold. Connected components of intersections of $\overline{\Cal{M}}_{0,S}(\bold{R})$
with complex boundary strata form a cell decomposition.
Cells of it are in one--to--one correspondence
with stable locally planar $S$--labeled trees.
The relation ``a cell is a codimension one
component of the boundary of another cell''
corresponds to the relation ``a locally planar tree
produces another locally planar tree by contracting
an internal edge.''

\smallskip

(ii) Fix an unoriented cyclic order on $S$ and consider the 
respective open cell. Any choice of three consecutive labels
with respect to this order allows one to introduce real
coordinates which
identify the open cell with the simplex $\Delta_n^0$,
where $|S|=n+3.$ This identification was denoted
$\Phi$ in (5).

\smallskip

(iii) The closure of each open cell
has the structure of a Stasheff polytope. In particular,
its boundary strata of codimension $1$ 
are indexed by those stable 2--partitions of $S$
which are compatible with the respective cyclic order:
they correspond to breaking the real equator into two connected arcs. 

\smallskip

(iv) The Zariski closure of the boundary of any open cell
is a union of boundary divisors having the same combinatorial shape
(``an algebraic Stasheff polytope'').}

\medskip

{\bf Sketch of a proof.} One can extract a proof of
this  Proposition from [K] (cf. also [Dev]). We will give here a modified
version including a differential geometric picture of the degeneration
of stable curves which  in particular makes it clear why
Stasheff polytopes appear at all. Assume that $|S|\ge 4.$

\smallskip

Let us look first at the open cells. A contemplation will
convince the reader that they are labeled by the cyclic orders $\rho$
of $S$ and that an arbitrary order  can arise. 
To identify each cell with  a simplex,
we choose three consecutive
labels with respect to $\rho$ and denote them
$1,\infty, 0$. For any $s\in S - \{1,\infty, 0\},$
consider the forgetful morphism 
$\overline{\Cal{M}}_{0,S}\to \overline{\Cal{M}}_{\{0,1,\infty ,s\}}.$
The target space is a projective line with
three special points: boundary divisors
corresponding to the partitions $s0|1\infty$, $s1|0\infty$, $s\infty |01$. 
\smallskip

We define $t_s$ as the lift from $\overline{\Cal{M}}_{0,\{0,1,\infty ,s\}}$
of the affine coordinate on this projective line taking respectively
the values $0,1,\infty$ at these three boundary divisors.

\smallskip

Put $S - \{1,\infty, 0\} =\{s_1,\dots ,s_n\}$ 
where the labels are cyclically ordered so that $0<s_1< \dots < s_n <1<\infty <0.$
From our definitions
it follows that with coordinates $t_{s_i}$, the interior part of the respective
real component can be identified with the simplex
$\Delta_n^0$ as in the Introduction.

\smallskip

\smallskip

Now let us turn to the boundary of the closure of this simplex
in $ \overline{\Cal{M}}_{0,S}(\bold{C}).$ As we have already said,
a point $z\in \Cal{M}_{0,S} (\bold{C})$
corresponds to the Riemannian
sphere $\bold{P}^1 (\bold{C})$ with marked points labeled by $S$.
If $z$ is real, this sphere has a distinguished equator $\bold{P}^1 (\bold{R})$
carrying all marked points $x_s(z)$.

\smallskip

Endow  $\bold{P}^1 (\bold{C}) - \{x_s(z)\,|\,s\in S\}$
with a hyperbolic metric of constant curvature $-1$.
The marked points will go to infinity as the ends
of infinitely thinning tubes, and the whole picture
will resemble a chestnut.

\smallskip

When $z$ tends to a point of a boundary divisorial stratum,
the respective variable surface carries a
closed geodesic whose length tends to
zero, a representative of the vanishing 1--homology class.
The body of the chestnut near the limit becomes divided into
two new bodies, connected by a thin tube whose length
tends to infinity as well. A part of labeled points will belong to
one half, the remaining ones to another half. This unordered
2--partition determines the relevant boundary component of
the moduli space.

\smallskip

The key fact is this: if $z$ tends to the boundary in this way
along {\it a real path} in $\Cal{M}_{0,S} (\bold{C})$, 
then the vanishing geodesic intersects
the real equator exactly at two points.  When the vanishing
cycle shrinks, the equator itself
degenerates into a union of two real equators
of the components of the limit curve. Hence  the respective 2--partition
of $S$ can be obtained
by breaking the equator into two arcs
(say, by drawing a chord).

\smallskip

From [K] it follows that all such stable 2--partitions 
compatible with the cyclic order arise in this way.
Geometrically, such a partition arises in the limit when all points
of one part of the partition try to merge.
 
\smallskip

Strata of the real boundary of larger codimension (faces)
can be obtained by iterating this description, i.~e.
applying it in turn to the components of a degenerate curve.
In this way, these strata become indexed 
by ``cactiform'' trees, whose ``edges'' are circles carrying marked points.
As in the complex case, we have applied dualization and replaced
 these graphs by
actual trees: vertices replace circles, halves of edges
replace marked points, edges replace contact points of pairs of
circles. 

\smallskip

Equivalently, put the points
labeled by $S$ on a circle in a plane compatibly with the chosen order.
Strata of codimension $d$ of the closure of the respective component
are in a natural bijection with the isotopy classes of stable trees
with $d$ (internal) edges
situated inside this circle and having the labeled
points as ends of tails. A stratum of codimension $d+1$
belongs to the closure of the stratum of codimension $d$
if one can obtain the smaller tree from the larger one by collapsing
one of its edges. The resulting polytope is called
the Stasheff polytope $K_n$ (if $|S|=n+3$). Its vertices correspond to the
trivalent trees.

\smallskip

Notice that in the Proposition 2.1 we were speaking about locally planar trees
whereas the picture above refers to a global planar embedding.
The point is that any given cell belongs to the boundaries of many
maximal cells. They are classified precisely by those planar embeddings
which are compatible with a given locally planar structure.

\medskip

\input epsf
\epsfxsize = 2in
\midinsert
$$\centerline{\hbox{\epsffile{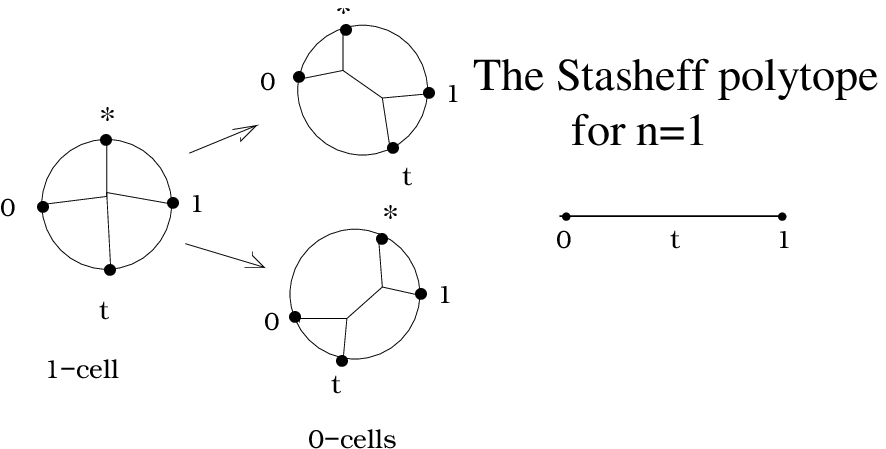}}}$$
\endinsert

\medskip
\epsfxsize = 2.5in
\midinsert
$$\centerline{\hbox{\epsffile{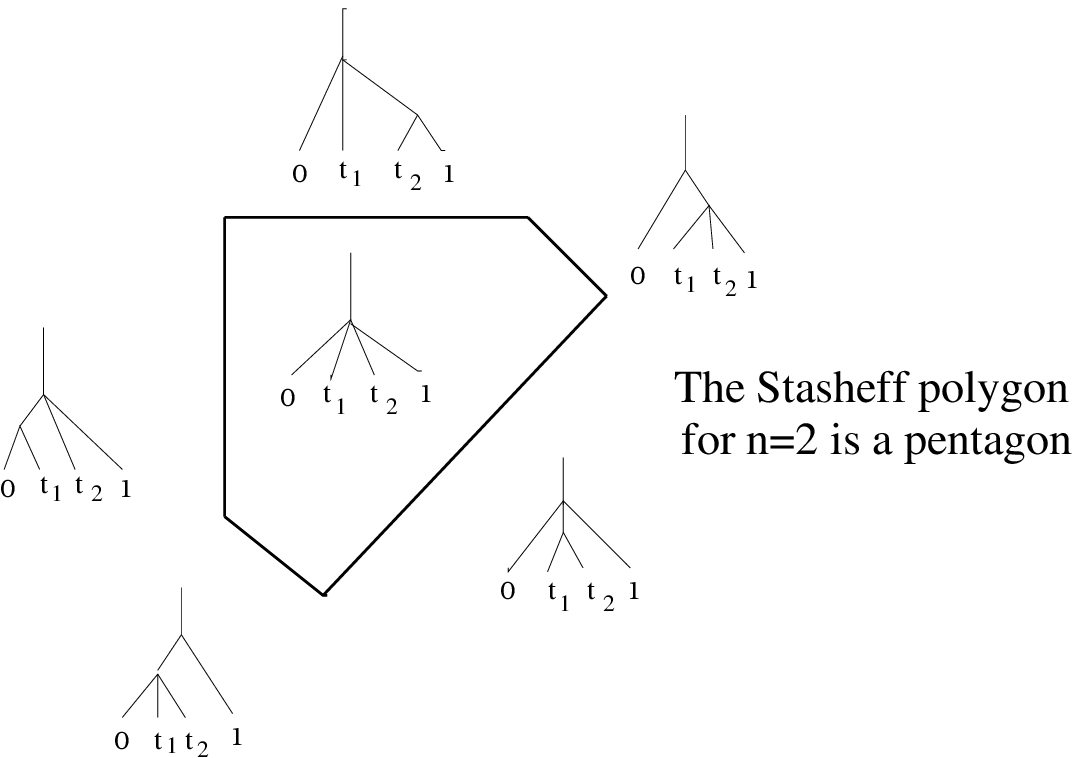}}}$$
\endinsert

\medskip

%\bigskip

%\centerline{[HERE GO THE PICTURES OF $K_1$ AND $K_2$]}

%\bigskip

Any face of  a Stasheff polytope is a product of Stasheff
polytopes (some of them may be points).

\smallskip

The Zariski closure $B$ of the boundary of any open cell 
is a union of boundary divisors having the same combinatorial shape
as a Stasheff polytope $K$ in the following sense: there is a bijection between the irreducible
components $D_i$ of $B$ and the codimension one faces $F_i$ of $K$
such that a subset of $B_i$'s has a non--empty intersection of expected codimension
iff the respective subset of $F_i$'s has this property.

\smallskip

A more elementary example: $n+1$ generic hyperplanes in
$\bold{P}^n$ have the shape of a simplex. However, four lines in
$\bold{P}^2$ do not have shape of a convex quadrangle.   

\medskip

{\bf 4. Divisors $B(\rho )$.} It is convenient to
repeat the formal definition of $B(\rho )$: if $\rho$
is any cyclic order on $S$, $B(\rho )$ is the sum
of boundary divisors corresponding to those 2--partitions
that are provided by breaking the circle into two arcs.
This divisor carries the total boundary of the
respective open cell $\Cal{M}^{\rho}_{0,S}(\bold{R}).$

\bigskip

{\bf 3\quad Forms  $\Omega (\underline \varepsilon)$ and divisors 
$A(\underline \varepsilon )$}

\medskip

{\bf 1. Construction of $\Omega (\underline \varepsilon )$.}
As above, choose and fix three labels in $S$ which we will
denote $0,1,\infty .$ We have shown that this choice allows us to define
a meromorphic function $t_s$ on $\overline{\Cal{M}}_{0,S}$
for each $s\in S - \{0,1,\infty \}.$ 

\smallskip

The additional structure is an arbitrary function
$\underline \varepsilon :\, S - \{0,1,\infty \} \to
\{0,1\}$,  $s\mapsto \varepsilon_s$. Choosing an
order of $S - \{0,1,\infty \}$, we can  define
a meromorphic form on $\overline{\Cal{M}}_{0,S}$:
$$
\Omega (\underline \varepsilon ):=
\bigwedge_{s\in S - \{0,1,\infty \}} \,\frac{dt_s}{t_s-\varepsilon_s}\, .
\eqno(7)
$$
Changing the order influences only the sign of the form.

\medskip

{\bf 2. Divisor $A (\underline \varepsilon )$.} We continue
using notation as in (7).

\smallskip

For $\alpha \in \{0,1,\infty\}$, we define $S(\alpha )= 
S(\alpha , \underline \varepsilon )$ by
$$
S(0):=\{s\,|\,\varepsilon_s=0\},\
S(1):=\{s\,|\,\varepsilon_s=1\},\
S(\infty ):= S(0)\,\cup\,S(1)\,.
$$

\smallskip

We will say that a 2--partition of $S$ {\it has type $\alpha$ with respect
to $\underline \varepsilon$} if one part
of it is of the form $\{\alpha\}\cup T$
where $T$ is a non--empty subset of $S(\alpha ).$

\medskip

{\bf Proposition 3.1.} {\it The divisor of singularities
 of $\Omega (\underline \varepsilon )$ 
on $\overline{\Cal{M}}_{0,S}$ is the sum
$A(\underline \varepsilon )$ of boundary
divisors corresponding
to those stable 2--partitions of $S$ which have some type with respect to
$\underline \varepsilon$.}

\smallskip

The following result is stronger than Lemma 1.1.
Denote by $\rho$ the cyclic order $0< s_1< \dots s_n <1 <\infty < 0.$

\medskip

{\bf Corollary 3.2.} {\it Assume that $n\ge 2$ and $\underline \varepsilon$
and $\rho$ satisfy $\varepsilon_{s_1}=1,  \varepsilon_{s_n}=0.$
Then $A (\underline \varepsilon )$ does not intersect the boundary of
$\Cal{M}^{\rho}_{0,S}(\bold{R}).$ 

\smallskip

In particular, $\Omega (\underline \varepsilon )$
is holomorphic on $\overline{\Cal{M}}_{0,S} - A(\underline \varepsilon )$, 
the closure of $\Phi (\Delta_n)$ determines a relative cycle on
$\overline{\Cal{M}}_{0,S} - A(\underline \varepsilon )$ modulo
$B(\rho )\,\cap\,A(\underline \varepsilon ) $, 
and formula (5) can be written as the integral over this relative cycle. }  

\medskip

{\bf Deduction of the Corollary.} Assume that the intersection in question
is not empty. Since it is closed, it has to contain
a zero--dimensional face of (the closure of) $\Cal{M}^{\rho}_{0,S}(\bold{R}).$
As was explained above, this face 
must correspond to a planar tree $\tau$ having all vertices
of multiplicity three and drawn inside the circle
on which the labels $S$ are put in the order $\rho$ and mark
the tails of $\tau$. The same tree describes the respective
complex zero--dimensional stratum of $\overline{\Cal{M}}_{0,S}.$

\smallskip 

Each boundary divisor containing such a stratum can be obtained by choosing
an interior edge of $\tau$: if we cut the edge, $S$
will break into two parts corresponding to the distribution
of tails among the resulting connected components of $\tau$.
Clearly, the same 2--partition of $S$ can be obtained
by breaking the circle into two connected arcs. 

\smallskip

It remains to check  that such a partition
cannot have a definite type with respect to $\underline{\varepsilon}.$
In other words,
no arc  can contain
only points of $S$ forming a set $\{\alpha \}\,\cup\, T$, 
$T\subset S(\alpha ).$
In fact, for $\alpha =\infty$ this is impossible,
because $1$ and $0$ separate $\infty$ from any element
$s_i$. For $\alpha = 0$ this is impossible because
$s_1\in S(1)$ and $\infty$ separate $0$ from any
element of $S(0).$ The case $\alpha =1$ is impossible
for similar reason.

\medskip

{\bf 3. Proof of the Proposition 3.1.} We will argue by induction on $n$.

\smallskip

{\it Case $n=1$.} Here $S=\{0,1,\infty, s\}.$ As we recalled in Chapter 2,
$\overline{\Cal{M}}_{0,S}$ is a projective line with
three  boundary divisors
corresponding to the partitions $s0|1\infty$, $s1|0\infty$, $s\infty |01$.
Moreover, $t:=t_s$ is the affine coordinate on this projective line taking respectively
the values $0,1,\infty$ at these three boundary divisors.
A form $\dfrac{dt}{t-\varepsilon},\, \varepsilon =0,1,$  is singular
at two partitions: $s\varepsilon|(1-\varepsilon)\infty$ and  $s\infty |01$.
Comparing this with the definition in the first subsection of this
Chapter, we see that these are exactly those partitions that have
a definite type with respect to $\varepsilon$.

\medskip

{\it Case $n=2$.} Although it is not strictly necessary to treat this case
separately, it will be convenient to do so:
we will show both the essence of the inductive step and 
the geometry of the first simple but nontrivial situation.

\smallskip

We put $S=\{0,1,\infty, s_1,s_2\}$, $t_i=t_{s_i}$, and choose 
$\varepsilon_1 =1, \varepsilon_0=0$. For $n=2$, this is the only case producing
a $\zeta$--value, namely, $\zeta (2)$.

\smallskip

Two morphisms forgetting $x_{s_2}$ and $x_{s_1}$ respectively,
represent $\overline{\Cal{M}}_{0,S}$
as a blow up of a quadric:
$$
\beta :\,\overline{\Cal{M}}_{0,S} \to \overline{\Cal{M}}_{0,\{0,1,\infty, s_1\}} \times
\overline{\Cal{M}}_{0,\{0,1,\infty, s_2\}} \cong \bold{P}^1\times \bold{P}^1\, .
\eqno(8)
$$
In coordinates $(t_1,t_2)$ three points are blown:
$(0,0),\,(1.1),$ and $(\infty ,\infty ).$ The resulting divisors
are (correspond to) respectively the following partitions of $S$: 
$s_1s_20|1\infty$, $s_1s_21|0\infty$, $s_1s_2\infty |01$.
The strict preimage of the diagonal is the boundary
divisor $s_1s_2|01\infty$. All in all, there are ten boundary components. We split them into a union of
two ``pentagons''. The pentagon $A=A(1,0)$ consists of the strict
preimages of the divisors
$t_1=1,t_1=\infty, t_2=0, t_2=\infty ,$ and of the blow up
of $(\infty ,\infty ).$ 

\smallskip

This is the picture we referred to
in the Introduction.

\smallskip

It is important to understand why the pentagon $A$ is precisely
the divisor of singularities of
$$
\Omega (1,0) = \beta^* \left( \dfrac{dt_1}{t_1-1}\,\wedge \,\dfrac{dt_2}{t_2}\right) .
$$
The point is that if we naively add up the total preimages of singularities
of $\dfrac{dt_1}{t_1-1}$ and $\dfrac{dt_2}{t_2}$, we will 
get two spurious extra components: blow ups of $(0,0)$ and
$(1,1)$.  A local calculation
shows that they cancel in $\Omega (1,0)$, because neither $(0,0)$ nor $(1,1)$
is a stratum of the divisor of singularities 
of the form on the quadric.
 
\smallskip

It is this step that is crucial -- and somewhat cumbersome to check --
in the general inductive reasoning.

\smallskip

To complete the discussion of our
example, notice that the complementary pentagon which is the union of the strict
preimages of $t_1=0,t_1=t_2,t_2=1$ and the blow ups
of $(0,0)$ and $(1,1)$, is precisely $B$. The reader may 
check it as a nice exercise. 

\medskip

{\it General inductive step.} Let $S=\{s_1,\dots ,s_n\}$,
$S'= \{s_1,\dots ,s_n, s_{n+1}\}$, $\varepsilon_{s_i}=\varepsilon_i$,
$\underline \varepsilon =\{\varepsilon_1,\dots ,\varepsilon_n\}$,
$\underline \varepsilon' =\{\varepsilon_1,\dots ,\varepsilon_{n+1}\}$.
Generalizing (8), we define a birational morphism
composed of two forgetful maps
$$
\beta :\,\overline{\Cal{M}}_{0,S'} \to \overline{\Cal{M}}_{0,S} \times
\overline{\Cal{M}}_{0,\{0,1,\infty, s_{n+1}\}} \, .
\eqno(9)
$$
From (7), it is clear that
$$
\Omega (\underline \varepsilon')=\beta^*\left(
 \Omega (\underline \varepsilon )\, \wedge\,\frac{dt_{n+1}}{t_{n+1}-\varepsilon_{n+1}}
\right)\, .
\eqno(10)
$$
The divisor of (logarithmic) singularities of $\Omega (\underline \varepsilon )\, \wedge\,\frac{dt_{n+1}}{t_{n+1}-\varepsilon_{n+1}}$ is a normal crossing
divisor consisting of (lifts of) some boundary divisors of the factors.
More precisely, by the inductive assumption, the latter constitute
$A(\varepsilon_{n+1})$ and $A(\underline \varepsilon )$ respectively.

\smallskip

Hence the divisor of $\Omega (\underline \varepsilon')$ is
$\beta^*(\roman{pr}_1^*A(\underline \varepsilon )+ \roman{pr}_2^*A(\varepsilon_{n+1}))$
minus eventual spurious components.

\smallskip

According to the Lemma 3.8 in [G5], the spurious components are
those irreducible boundary divisors in $\overline{\Cal{M}}_{0,S'}$
which get blown down by $\beta$ to a subvariety
of the product {\it which is not a stratum of the
divisor} $\roman{pr}_1^*A(\underline \varepsilon )+ \roman{pr}_2^*A(\varepsilon_{n+1}).$
Let us now make explicit the combinatorics of the situation.
We will systematically identify boundary divisors
with partitions.

\smallskip

(i) First of all, a boundary divisor of $\overline{\Cal{M}}_{0,S'}$ gets 
blown down by $\beta$ if both forgetful morphisms map it into 
boundary divisors. Combinatorially, we take a stable partition of
$S'$ and induce from it a partition of $S$ and 
of $\{0,1,\infty ,s_{n+1}\}$; both results must be stable.
A contemplation will convince the reader that
there are three kinds of such partitions:
$$
s_{n+1}0\sigma_1|\sigma_21\infty ,\ 
s_{n+1}1\sigma_1|\sigma_20\infty ,\
s_{n+1}\infty\sigma_1|\sigma_201\,
\eqno(11)
$$
where $\sigma_1|\sigma_2$ is any partition of $\{s_1,\dots ,s_n\}$
with $\sigma_1\ne \emptyset .$

\smallskip

Moreover, divisors $s_{n+1}\sigma_1|\sigma_201\infty$ 
with $|\sigma_1|\ge 2$ also
get blown down, and this completes the list: see [Ke],
Lemma 1, p. 554.

\medskip

(ii) Now let us list components of $A(\underline \varepsilon ).$
They are also divided into three types (see the beginning of the Chapter):
$$
0T_0|\dots 1\infty ,\  1T_1|\dots 0\infty ,\  \infty T_{\infty}|\dots 01
\eqno(12)
$$
where $T_{\alpha}$ is a non--empty subset of $S(\alpha ).$
Each of them has exactly two strict lifts to $\overline{\Cal{M}}_{0,S'}$:
we may add $s_{n+1}$ to the first or to the second
parts of the partition. We will first check which
of these lifts {\it do not} belong to $A(\underline \varepsilon').$
Assume for concreteness that $\varepsilon_{n+1}=0$; the other case
is treated similarly.
\smallskip

Again, a contemplation will convince the reader
that only the following lifts do not have a definite type
with respect to $\underline \varepsilon'$:
$$
s_{n+1}1T_1|\dots 0\infty \, .
\eqno(13)
$$

\smallskip

(iii) Let us now check that the components (13) are spurious.
First, forgetting $s_{n+1}$ or all $\{s_1,\dots ,s_n\}$
produces stable partitions, so such a divisor
gets blown down. More precisely, its $\beta$--image is
(not just lies in)
$$
1T_1|\dots 0\infty \times s_{n+1}1| 0\infty
$$
This product is not a stratum of the divisor
of $\Omega (\underline \varepsilon )\,\wedge\, \dfrac{dt_{n+1}}{t^{n+1}}\,.$ 
In fact, the first factor is a component of 
the singularities of $\Omega (\underline \varepsilon )$,
whereas the second one is not a singularity for
$\dfrac{dt_{n+1}}{t_{n+1}}$.

\smallskip

(iv) After adding $s_{n+1}$ to one of the parts of partitions
in the list (12) and deleting the partitions (13) we will
get almost all partitions of $S'$ having a definite type
with respect to $\underline \varepsilon'$.
The only exceptions will be
(recall that $\varepsilon_{n+1}=0$)
$$
s_{n+1}0\,|\, S1\infty , \quad s_{n+1}\infty\,|\, S01\, .
$$
These two components of $A(\underline \varepsilon')$
are supplied by the lifts of the poles of
$\dfrac{dt_{n+1}}{t_{n+1}}$.

\smallskip

This concludes the proof.

\bigskip

{\bf 4\quad Multiple $\zeta$--motives}

\medskip

{\bf 1. Motivic background.}  Below we employ the abelian category ${\Cal
M}_T(\bold{Q})$
of mixed Tate motives over $\bold{Q}$ as defined in chapter 5 of
[G3], (see also [L1]), and its subcategory ${\Cal M}_T(\bold{Z})$ of
unramified  over $\bold{Z}$
mixed
Tate motives  ([G4], [DG]).
As is shown in loc.~cit., an
object of ${\Cal M}_T(\bold{Q})$ is unramified over $\bold{Z}$ if and only if for
any prime
$l$ its $l$--adic realization is unramified outside of $l$.

\smallskip

To talk about the periods of mixed Tate motives
we need to equip mixed Tate motives with an additional structure called
framing.
For a definition of  framed mixed Tate objects, e.~g.
framed mixed Tate motives or framed Hodge--Tate structures,
see for example section 3.2 of [G4]. To define a  period of a
framed Hodge--Tate structure $H$ we need to choose
a splitting of the weight filtration on $H_{\bold{Q}}$, see 
the next section and section 4.2 in
[G3] for more details.

\medskip

{\bf 2. Multiple $\zeta$--motives.} Keeping the notation of section 1, Chapter 2, we
define  
$$
\zeta^{\Cal{M}}(n_1,\dots ,n_m) :=
H^n(\overline{\Cal{M}}_{0,n+3} - A(\varepsilon ), B - B\,\cap\,A(\varepsilon ) ) .
\eqno(14)
$$

\medskip

{\bf Theorem 4.1.} {\it (i) $\zeta^{\Cal M}(n_1, ..., n_m)$ is
a mixed Tate motive unramified over $\roman{Spec}\,(\bold{Z})$.

\smallskip

(ii) $\Omega (\underline \varepsilon )$ and $\Phi (\rho )$ allow us to 
define its natural framing so that the respective  period
equals $\zeta(n_1, ..., n_m)$.}

\smallskip

{\bf Proof.} A non singular projective variety is called a Tate variety if its motive
is a direct
sum  of  pure Tate motives.  It is well known (and follows 
from Knudsen's inductive
description of ${\overline {\Cal M}}_{0, n+3}$) 
that all strata of the stratification
defined by the
 divisor $\partial {\overline{\Cal M}}_{0, n+3}$ are  Tate
varieties.
Since
$\partial {\overline{\Cal M}}_{0, n+3}$ is a normal crossing divisor,
Corollary 3.2 above and Proposition 3.6 in [G5] show that 
$\zeta^{\Cal M}(n_1, ..., n_m)$ is well defined as  a mixed Tate motive
over $\bold{Q}$.

\smallskip

Let us recall briefly its construction.
Consider the standard cosimplicial variety
$$
S_{\bullet}({\overline{\Cal M}}_{0,n+3} -A(\underline \varepsilon), B)
\eqno(15)
$$ 
corresponding to the pair
$(\overline{\Cal{M}}_{0,n+3} -A(\underline \varepsilon),
B-B\,\cap\,A(\underline \varepsilon))$. 

\smallskip

Recall that
$S_{0} := \overline{\Cal{M}}_{0,n+3} -A(\underline \varepsilon)$, and  $S_{k}$ is the
disjoint union
of the codimension $k$ strata of the divisor $B-B\,\cap\,{A(\underline \varepsilon)}$.

\smallskip

Following the standard procedure, let us make a complex of varieties out of
(15) with
$S_0$ placed at the degree $0$.
It provides an object in the Voevodsky triangulated category of motives [V],
which in fact belongs to the
triangulated subcategory ${\Cal D}_T(\bold{Q})$ of mixed Tate motives over $\bold{Q}$.
There exists a
canonical $t$--structure $t$ on the category ${\Cal D}_T(\bold{Q})$ (see Chapter 5 of
[G3] or [L1] for its definition).
Then $H^n_t$ of the above complex is the multiple $\zeta$--motive.

\smallskip

Now we will introduce a framing. Generally, a framing on a mixed
motive is an algebraic
counterpart of the notion of  period. For the mixed Tate motives it
looks as follows.
Recall that a mixed Tate motive carries a canonical weight filtration
$W_{\bullet}$ such that the associated graded motive of weight $-2k$ is a direct
sum of
pure Tate motives $\bold{Q}(k)$, and ${\roman Gr}^W_{2k-1} =0$. An $n$--framed (or
simply framed if $n$ is clear from the context) mixed Tate motive is a triple $(M, v,
f)$ where $v$ and $f$ are two non zero morphisms
$$
v: \,\bold{Q}(-n) \to {\roman Gr}^W_{2n}M\, , \quad f:\, \bold{Q}(0)\to
\left({\roman Gr}_{0}^WM\right)^{\vee} = {\roman Gr}_{0}^W M^{\vee}
$$
and $M^{\vee}$ is the dual object.

\smallskip

In our situation we  define a framing as follows.
There are canonical isomorphisms
$$ 
{\roman Gr}^{W}_{2n}H^n(\overline {\Cal M}_{0,n+3} -A(\underline
\varepsilon),
B -B\,\cap\,A(\varepsilon )) \cong
{\roman Gr}^{W}_{2n}H^n(\overline {\Cal M}_{0,n+3} -A(\underline
\varepsilon))
\eqno(16)
$$
and
$$ 
{\roman Gr}^{W}_{0}H_n(\overline {\Cal M}_{0,n+3} -A(\underline
\varepsilon), B -B\,\cap\,A(\varepsilon ) \cong
{\roman Gr}^{W}_{0}H_n(\overline {\Cal M}_{0,n+3}, B)
\eqno(17)
$$
Moreover, there are natural non--zero morphisms of pure Tate motives
$$
[\Omega(\underline \varepsilon)]:\, \bold{Q}(-n) \to {\roman
Gr}^{W}_{2n}H^n(\overline {\Cal M}_{0,n+3} -A(\underline \varepsilon))\, ,
$$
$$
[\Phi_n]: \bold{Q}(0) \,\to {\roman Gr}^{W}_{0}H_n(\overline {\Cal M}_{0,n+3}, B) \,.
$$
To construct them we employ the fact that the Hodge realization is a
fully faithful functor
on the category of pure Tate motives. Then the first morphism is
determined
in the De Rham realization by the form
$\Omega(\underline \varepsilon)$. The second
morphism is determined in the Betti realization by the relative homology
class of the cell
$\overline \Phi_n$.  Combining them
with the isomorphisms (16) and (17) we get the
frame morphisms
$$
[\Omega(\underline \varepsilon)]': \,\bold{Q}(-n) \to {\roman
Gr}^{W}_{2n}H^n(\overline {\Cal M}_{0,n+3} -A(\underline \varepsilon),
B-B\,\cap\,A(\varepsilon ))\, ,
$$
$$
[\Phi_n]':\, \bold{Q}(0) \to {\roman Gr}^{W}_{0}H_n(\overline {\Cal M}_{0,n+3}
-A(\underline
\varepsilon), B-B\,\cap\,A(\varepsilon ))\, .
$$
Let us show that $\zeta(n_1, ..., n_m)$ is a period
of this framed mixed Tate motive.

\smallskip

Recall that a mixed Hodge structure is the following linear algebra
data:

i) A $\bold{Q}$--vector space $H_{\bold{Q}}$ equipped with   an increasing
filtration $W_{\bullet}$.

\smallskip

ii) A decreasing filtration $F^{\bullet}$ on $H_{\bold{C}} = H_{\bold{Q}}\otimes \bold{C}$
inducing for each integer $n$  a weight $n$ pure Hodge structure
on ${\roman Gr}^W_nH$, i.e.
$$
{\roman Gr}^W_nH_{\bold{C}} = \oplus_{p+q=n}F^p{\roman Gr}^W_nH_{\bold{C}} \cap \overline
F^q{\roman Gr}^W_nH_{\bold{C}}
$$
A Hodge--Tate structure is a mixed Hodge structure with the Hodge numbers
$h^{pq}=0$ unless $p=q$. It is easy to see that this is equivalent
to the following condition: for every $p\in \bold{Z}$ the natural map
$$ 
F^pH_{\bold{C}} \cap W_{2p}H_{\bold{C}} \to {\roman Gr}_{2p}^WH_{\bold{C}}
\eqno(18)
$$
is an
isomorphism.

\smallskip

The Hodge realization of the multiple $\zeta$--motive is a Hodge--Tate
structure.

\smallskip

An $n$--framed
Hodge--Tate structure is a Hodge--Tate structure $H$ equipped with
non--zero morphisms
$$
v: \bold{Q}(-n) \to {\roman Gr}^W_{2n}H; \qquad f: \bold{Q}(0)\to
\left({\roman Gr}_{0}^WH\right)^{\vee} = {\roman Gr}_{0}^W H^{\vee}
$$

To define a period of a framed Hodge-Tate structure we need
to choose an additional data -- a map of $\bold{Q}$--vector spaces
$\widetilde f:\, \bold{Q} \to H_{\bold{Q}}^{\vee}$ lifting  $f$, i.e.
${\roman Gr}_{0}^W \widetilde f = f$.
Set
$f':= \widetilde f(1)$.
The composition
$$
\bold{Q}(-n) \to {\roman Gr}_{2p}^WH_{\bold{Q}} \to F^pH_{\bold{C}} \cap
W_{2p}H_{\bold{C}}
$$
where the first arrow is $v$ and the second arrow is provided by  (18), leads to
a vector $v' \in F^pH_{\bold{C}} \cap W_{2p}H_{\bold{C}}$.
The corresponding period
 is a complex number $\langle v',  f'\rangle$. A different choice
of lifting $\widetilde f$  changes the period by $2 \pi i$ times a
``weight $n-1$ period'', see
s. 4.2 in [G3] for more details.

\smallskip

In our case  a canonical choice of the
lift $\widetilde f$ is secured by the Corollary 3.2
which shows that
the cell $\overline \Phi_n$ provides an element of the Betti
homology
$$
[\overline \Phi_n]' \in H^B_n(\overline {\Cal M}_{0,n+3} -A(\underline
\varepsilon),
B - B \cap {A(\underline \varepsilon)})\, .
$$
The restriction of the form $\Omega(\underline \varepsilon)$ to the
divisor $B$ is zero, so that it
furnishes  a  De Rham cohomology class
$$
[\Omega(\underline \varepsilon)]'\in H_{\roman DR}^n(\overline {\Cal
M}_{0,n+3} -A(\underline \varepsilon),
B - B \cap {A(\underline \varepsilon)})\, .
$$
These are the classes obtained before  by using isomorphisms
(16) and (17).
Clearly, the integral (5) computes the pairing of these two
classes.

\smallskip

It remains to show that $\zeta^{\Cal M}(n_1, ..., n_m)$ is unramified
over $\bold{Z}$.

\medskip

{\bf Definition 4.2.} {\it Let $D$ be a normal crossing divisor
in a regular scheme $X$ over $\bold{Z}_p$. Assume that the pair
$(D,X)$ is proper over $\bold{Z}_p$.
We say that reduction modulo $p$ does not change the combinatorics of
$(D,X)$ if
$X$ and every stratum of $D$ are smooth over $\bold{Z}_p$, and
the reduction map from the strata of $D$ to ones at the special fiber
is a bijection.}

\smallskip

%Let $Y$ be an integral scheme  of finite type over $\bold{Z}_p$ whose generic
%fiber
%$Y_{\eta}$ is non empty. Let $Y_0$ be
%its special fiber.  Then $Y$ is smooth over $\bold{Z}_p$
%if and only if $Y_0(\overline F_p)$ and $Y_{\eta}(\overline{\bold{Q}}_p)$
%are nonsingular
%(see Proposition 2.9 in the Chapter IV of [S]).

\smallskip
Let $X$ be a flat integral scheme  of finite type over $\bold{Z}_p$, $D$ a divisor.
Set $\overline X:= X \otimes_{\bold{Z}_p} \overline{\bold{Q}}_p$.
Let $A$ and $B$ be the unions
of two disjoint subsets of irreducible components of $D$.

\medskip

{\bf Proposition 4.3.} {\it If the reduction
modulo $p$ does not change the combinatorics of $(D,X)$ and $(l,p)=1$,
then the
$\roman{Gal}(\overline{\bold{Q}}_p/\bold{Q}_p)$--module
$
H^n_{\roman{et}}(\overline X - \overline A, \overline{B}-
\overline{B}\,\cap\,\overline{A}; \bold{Q}_l)
$ is unramified at $p$.}

\smallskip

{\bf Proof}. It is proved by taking the standard simplicial resolution
computing $H^n(X-A,B- B\,\cap\,A)$ and using the proper and  smooth base change
theorem,
see Proposition 3.9 in [G5] for details.

\medskip

{\bf Proposition 4.4.} {\it For any prime $p$ the reduction
modulo $p$ does not change the combinatorics of the divisor
$\partial \overline {\Cal M}_{0,n+3}$.}

\smallskip

In fact, the description of the combinatorics of the boundary
strata given in the subsection 2.2 is valid over a field
of any characteristic $p$ and is compatible
with reduction. 

\smallskip

This concludes the proof of the Theorem 4.1.

\medskip

{\bf 3. Concluding remarks.} The discussion above can
be applied to more general pairs $A$, $B$ sharing
no common components for which we can find framing
as in (3). For example, interchanging $A$ and $B$ and using
the duality isomorphism
$$
H^n(\overline {\Cal M}_{0, n+3}-A)^{\vee} \cong
H_n(\overline {\Cal M}_{0, n+3}, A)(n)
$$
which implies
$$
\roman{Gr}^W_{2n}H^n(\overline {\Cal M}_{0, n+3}-A) \cong
\roman{Gr}^W_{0}H_n(\overline {\Cal M}_{0, n+3}, A)(n)
$$
we get a different
mixed Tate motive over $\bold{Z}$ with the dual framing.
Hence we get a concrete

\medskip

{\bf Problem}. {\it Let $\widehat \zeta^{\Cal M}(n_1, ..., n_m)$
be the
mixed Tate motive  obtained by interchanging the $A$- and $B$-divisors
in (14).  Express $\widehat \zeta^{\Cal M}(n_1, ...,
n_m)$ via the multiple $\zeta$--motives.}

\smallskip

%The equivalence classes of  unramified over $\bold{Z}$ framed mixed Tate
%motives
%form a graded Hopf
%algebra ${\Cal A}_{\bullet}(\bold{Z})$ (see section  3.2 of [G4]). The hypothetical
%answer
%to the above question  involves the coproduct of the
%element $\zeta^{\Cal M}(n_1, ..., n_m) \in \Cal{A}_n(\bold{Z})$.

\smallskip

According to the conjecture 17 b) in [G1]
any unramified over $\bold{Z}$ framed mixed Tate motive is equivalent to a
linear combination of
 framed multiple $\zeta$--motives. 
In our situation, we come to the following

\medskip

{\bf Conjecture 4.5.} {\it Any mixed Tate motive (2)
with a framing provided by (3)
is equivalent to a $\bold{Q}$--linear combination of the
framed weight $n$ multiple $\zeta$-motives.}

\smallskip

The analytic version of this conjecture says that for any
$$
[\Delta'_B] \in H_n(\overline {\Cal M}_{0, n+3}- A, B-B\,\cap\,A) \quad
\roman{lifting}
\quad [\Delta_B]\in \roman{Gr}^W_0H_n(\overline {\Cal M}_{0, n+3}, B)
$$
in (3)
the integral $(2 \pi i)^{-n}\int_{\Delta'_B}\omega_A$ giving a period
of the framed Hodge-Tate structure (4) is $(2 \pi
i)^{-n}\times$
(a $\bold{Q}$--linear combination
of the weight $n$ multiple $\zeta$--values).

\medskip

{\bf 4. Generalizations}. Consider an iterated integral
$$
\roman{I}(a_1, \dots , a_n):= \int_{\Delta_n^0}\frac{dt_1}{t_1-a_1} \wedge
\dots
 \wedge \frac{dt_n}{t_n-a_n}
\qquad a_1 \not = 0,\, a_n \not = 1\, .
\eqno(19)
$$
Let $\underline a:= \{a_1, \dots , a_n\}$ and $A(\underline a)$ be the
divisor of singularities
of the form
$$
\Omega(\underline a):= \frac{dt_1}{t_1-a_1} \wedge \dots \wedge
\frac{dt_n}{t_n-a_n}
$$
meromorphically extended to $\overline {\Cal M}_{0,n+3}$. Then
extending the definitions above in an obvious way,  one can
prove the following

\medskip

{\bf Proposition 4.6.} {\it 
$A(\underline a) \cup B$
is a normal crossing divisor and $A(\underline a)$ and $B$ share no
common components.}

\smallskip

The absence of common
components follows immediately from the convergence of the integral
(19).

\smallskip

Thus if $a_i$ are elements of a number field $F$ then
$$
H^n(\overline {\Cal M}_{0,n+3} - A(\underline a), B- B\,\cap\,{A(\underline a)})
$$
is  a mixed Tate motive over $F$.
The class $[\Delta_B]$ and the form $\Omega(\underline a)$ provide it
with a framing.
If $a_i$ are complex numbers then the period of the
corresponding framed Hodge--Tate structure is given by the integral
$\roman{I}(a_1, \dots , a_n)$.

\smallskip

This provides a different construction of the motivic multiple
polylogarithms motives
then the one suggested in [G5]. However
the framed mixed Tate motives given by these two constructions are
equivalent.

%\bigskip

\newpage

{\bf REFERENCES}

\medskip

 [DGr] Deligne P. {\it Resum{\'e} des premiers
expos{\'e}s de A. Grothendieck}.
SGA 7, LNM 288, 1--24.

\smallskip

 [DG] Deligne P., Goncharov A.B.  {\it (in preparation).}

\smallskip

[Dev] S.~L.~Devadoss. {\it Tesselations of moduli spaces
and the mosaic operad.} In: Contemp.~Math., vol. 239 (1999),
pp. 91--114.

\smallskip

 [Dr] Drinfeld V.G. {\it On quasitriangular quasi-Hopf algebras
and on a group that is closely connected with $\roman{Gal}
(\overline{\bold{Q}}/\bold{Q})$}. Leningrad Math. J. 2 (1991), no. 4, 829--860.

\smallskip

[E] L. Euler. {\it Opera Omnia.} Ser. 1, Vol XV, Teubner,
Berlin 1917, 217--267.

\smallskip
 
[G1] Goncharov A.B. {\it Polylogarithms in arithmetic and
geometry}.
Proc. of the International Congress of Mathematicians, Vol. 1, 2
(Z\"urich, 1994), 374--387, Birkh\"auser, Basel, 1995.

\smallskip

[G2] Goncharov A.B. {\it Multiple $\zeta$-values, Galois
groups and geometry of
modular varieties}. Proc. of the Third European Congress of
Mathematicians. Progress in Mathematics,
vol. 201, p. 361--392. Birkh\"auser Verlag, 2001.
math.AG/0005069.

\smallskip

[G3] Goncharov A.B. {\it Volumes of hyperbolic manifolds
and mixed Tate motives}. J. Amer. Math. Soc. 12 (1999) N2, 569--618.
math.alg-geom/9601021.

\smallskip

[G4] Goncharov A.B. {\it Multiple polylogarithms and mixed
    Tate motives}. \newline math.AG/0103059.

\smallskip

[G5] Goncharov A.B. {\it Periods and mixed motives}.
math.AG/0202154.

\smallskip

[K] Kapranov M.M. {\it The permutoassociahedron, MacLane's
    coherence theorem and
asymptotic zones for the KZ equation}. J. of Pure and Appl.
Algebra 85 (1993), 119--142.

\smallskip

[Ke] Keel S. {\it Intersection theory of moduli space
of stable $N$--pointed curves of genus zero}.
Trans.~AMS, 330:2 (1992), 545--574.

\smallskip

[Kn] Knudsen F.F. {\it The projectivity of the moduli space
    of stable curves II.
The stacks $\overline M_{0,n}$}. Math. Scand. 52 (1983), 163--199.

\smallskip

[Man]  Manin Yu.I. {\it Frobenius manifolds, quantum
cohomology, and moduli
spaces}. AMS Colloquium Publications, vol. 47, Providence, RI, 1999, 303
pp.

\smallskip

[L1] Levine, M. {\it Tate motives and the vanishing conjectures
for algebraic
$K$--theory}. Algebraic $K$--theory and algebraic topology (Lake
Louise, AB, 1991), 167--188, NATO Adv. Sci. Inst. Ser. C
Math. Phys. Sci., 407,
Kluwer Acad. Publ., Dordrecht, 1993.

\smallskip

[S] Silverman J. {\it Advanced topics in the arithmetic of
elliptic curves}.
Graduate texts in Mathematics, Springer, 1994.

\smallskip

[T] Terasoma T. {\it Multiple zeta values and mixed Tate
    motives}. math.AG/0104231.

\smallskip

[V] Voevodsky V. {\it Triangulated category of motives over a
field}. In:
Cycles, transfers, and motivic homology theories, 188--238, Ann.
of Math. Stud., 143, Princeton Univ. Press, Princeton, NJ, 2000.

\enddocument